\documentclass[12pt]{article}
\usepackage{amssymb}
\usepackage{amsmath}
\usepackage{graphicx}
\usepackage{latexsym}
\usepackage{tikz}
\def\part#1{\frac{\partial\phantom{q}}{\partial#1}}

\newenvironment{rmk}{\begin{trivlist}\item[]{\bf Remark:} }
{\end{trivlist}}

\newenvironment{ex}{\begin{trivlist}\item[]{\bf Example:} }
{\end{trivlist}}
\newenvironment{prf}{\begin{trivlist}\item[]{\bf Proof:} }
{\hfill $\Box$ \end{trivlist}}

\newtheorem{thm}{Theorem}

\newcommand{\lie}[1]{\mathfrak{#1}}

\def\OO{\mathcal O} 

\newcommand{\R}{\mathbf{R}}
\newcommand{\C}{\mathbf{C}}

\newcommand{\Z}{\mathbf{Z}}

\newcommand{\PP}{{\rm P}}

\begin{document}
\title{ ALE spaces and nodal curves}
\author{Nigel Hitchin\\[5pt]}
\maketitle
\section{Introduction} In \cite{Kron1} \cite{Kron2} Kronheimer gave a construction of  asymptotically locally Euclidean (ALE) metrics in four dimensions by   hyperk\"ahler quotients. These metrics are defined on the complex surfaces obtained by resolving the singularity at the origin of $\C^2/\Gamma$ where $\Gamma\subset SU(2)$ is a finite subgroup and the resolution is by a configuration of rational curves which intersect according to a Dynkin diagram of type ${\bf A,D,E}$. 
The construction  gives a great deal of global information about the metric -- completeness, asymptotic behaviour, moduli. Nevertheless, some aspects are impenetrable -- for example the case of ${\bf E_8}$ entails taking a vector space of dimension   $4\vert\Gamma\vert=480$ and solving $357$ quadratic equations even before taking a quotient to obtain a 4-manifold. 

The aim of this paper is to give a different approach based on twistor theory,  using algebraic geometry to highlight different ingredients, providing a fresh look  at these  spaces. Twistor theory was already used by the author 45 years ago  in the ${\bf A_k}$ and ${\bf D_k}$ cases  \cite{Hit1},\cite{CH} (even though the latter paper appeared much later) by using the algebraic equations defining the deformations of $\C^2/\Gamma$ embedded in $\C^3$ as a singular affine surface, or equivalently the other complex structures of the hyperk\"ahler family. This method produced some local formulae but proved impossible for  the consideration of the ${\bf E_k}$ series.

Here we replace the  equations by  a compactification of the twistor space. 
This aspect appeared in a recent paper of the author \cite{Hit3} and in some respects this one can be considered as a sequel. 
In both cases we consider the action of a circle  induced by  the scalar action of  $\C^*$ on $\C^2$ and consider the family of invariant metrics. The basis of twistor theory consists of studying a distinguished family of rational curves -- twistor lines -- in the twistor space.  In  \cite{Hit3} we identified the twistor lines which are preserved by the $\C^*$-action  in order to describe the induced metric on the central 2-sphere of fixed points in the ALE space. They are closures of $\C^*$-orbits. In this paper we identify the orbits  of a \emph{general} twistor line.  

The  twistor space is a complex 3-manifold which is a fibration $Z\rightarrow \PP^1$. The circle action  fixes two distinguished points $u=0,\infty$ on $\PP^1$ and the fibre over these points is isomorphic to  the complex structure of the resolved singularity. We introduce a compactification $\bar Z$ which extends the fibration and consider the fibre over $u=1$. This is a projective rational surface obtained by a systematic blowing up of points in the plane.  An orbit of a twistor line is a surface in $\bar Z$ which meets this fibre in a rational curve, and we identify in Theorem \ref{Q} the linear system to which these curves belong. The curves are however almost always singular,  so earlier  issues are transformed into the serious problem of identifying nodal curves in a known linear system. Nevertheless, we consider some examples -- the ${\bf A_{2\ell-1}}$ case involves hyperelliptic curves, and we also address the constraint on an elliptic curve which appeared  in the earlier treatment of the ${\bf D_4}$ case. 

It is at this point that we make contact with the work of N.Honda \cite{Honda1},\cite{Honda2} who adapted the original Penrose twistor approach by replacing smooth rational  curves with nodal ones. The relevant case here  consists of the 3-parameter family of  nodal curves in the  surface $u=1$ and this gives it an interpretation as the minitwistor space for an Einstein-Weyl structure 
 on the 3-dimensional quotient of the ALE space by the circle action.   

In the case of ${\bf A_1}$ -- the original Eguchi-Hanson metric which began the study of the locally Euclidean property -- the rational surface is a quadric and the curves have  no singularities -- they are just plane sections. As a consequence we deduce that the quotient of the ALE manifold can be identified with hyperbolic 3-space, and we use this as a model for a brief discussion of the orbifold Einstein-Weyl geometry in general. In particular, the fixed 2-sphere in the ALE space becomes a 2-sphere at  infinity in the quotient. Infinity in the ALE space becomes a distinguished point with the orbifold structure of $\R^3/\Gamma$. We do not develop further the properties of these spaces but the interpretation of a minitwistor space as the space of geodesics in an Einstein-Weyl manifold gives an explanation, in terms of closed geodesics, of the existence of the singularities in the rational  curves.

This paper is based on a talk given in Oxford on July 25th 2023 at a conference for Peter Kronheimer's 60th birthday. 
\section{The twistor construction}\label{twistor}
The ALE metrics are hyperk\"ahler, so have complex structures $I,J,K$ satisfying  the quaternionic identities and preserved by the Levi-Civita connection. They define  a 2-sphere of complex structures $aI+bJ+cK$ with $a^2+b^2+c^2=1$. To obtain a 4-dimensional hyperk\"ahler manifold via the twistor construction one needs a complex 3-manifold $Z$ with a projection $\pi:Z\rightarrow \PP^1$ and a family of sections which we call \emph{twistor lines} with normal bundle $N\cong \OO(1)\oplus \OO(1)$. Deformation theory gives a complex 4-dimensional space of such sections where the tangent space at a point is isomorphic to $H^0(P^1, N)$. A conformal structure is defined by declaring a null tangent vector to be a section of $N$ which vanishes at some point. In addition one requires a  symplectic form along the fibres -- a section of $\Lambda^2 T^*_{Z/\PP^1}(2)$ where $\OO(k)$ is the pull back by $\pi$ of the $\OO(k)$ on $\PP^1$. This fixes a complex metric in the conformal class. Finally we require a real structure, an antiholomorphic involution $\sigma:Z\rightarrow Z$ which covers the antipodal map $u\mapsto -1/\bar u$ on $\PP^1$, preserves the relative symplectic form and induces a positive definite metric on the space $M$ of real sections. Each fibre $Z_u$ meets a real section in a  unique point and endows $M$ with a complex structure and a complex symplectic form, with the fibre over $u\in \C\cup\infty= \PP^1\cong S^2$ realizing the 2-sphere of complex structures of the hyperk\"ahler family.  

Kronheimer's construction identifies the twistor space, confirming the picture suggested in \cite{Hit1}. A finite subgroup $\Gamma\subset SU(2)$ acts on $\C^2$ and there is a basis $x,y,z$ of the ring of homogeneous invariant polynomials in $(z_1,z_2)$ of degree $\ell,m,n$  respectively satisfying an equation $f(x,y,z)=0$ which  embeds the quotient $\C^2/\Gamma$ as a singular affine surface in $\C^3$. A versal deformation consists of adding lower degree polynomials to $f$, the coefficients of which can be understood as Weyl-invariant polynomials on the Cartan subalgebra $\lie{h}$ of the Lie algebra of type $ADE$ associated to $\Gamma$. Then  $f(x,y,z,h_1,\dots,h_k)=0$ is a family of affine surfaces, the generic one smooth. 
The twistor space for an $S^1$-invariant metric is the hypersurface in the total space of $\OO(\ell)\oplus \OO(m)\oplus \OO(n)$ defined by $f(x,y,z,a_1u,\dots,a_ku)=0$ with $a_i$ real, together with a resolution of the singularities at $u=0,\infty$.  The real structure is induced by the involution $(z_1,z_2)\mapsto (\bar z_2,-\bar z_1)$ on the invariant polynomials $x,y,z$ and $u=z_1/z_2$. 

The ${\bf A}$ and ${\bf D}$ series are easy to write down because the invariant polynomials are so well known, in topology for example as Chern classes for ${\bf A_k}$ and Pontryagin and Euler classes for ${\bf D_k}$:
\begin{eqnarray}
{\bf A_k}&:& xy=\prod_{i=1}^{k+1}(z-a_iu) \label{Aequation}\\
{\bf D_k}&:& x^2-zy^2=-\frac{1}{z}\left(\prod_{i=1}^k(z+a_i^2u^2)-u^{2k}\prod _{i=1}^k a_i^2\right)+2yu^k\prod_{i=1}^ka_i.  \label{Dequation}
\end{eqnarray}
Then finding the twistor lines is the ``Diophantine" problem of solving an algebraic equation with polynomial coefficients in $u$ with polynomials $x(u),y(u),z(u)$. 

The first case is straightforward and was carried out in \cite{Hit1}: $z(u)=cu^2+au-\bar c$  is a real quadratic  and  $z-a_iu=c(u-\alpha_i)(u-\beta_i)$. Since the discriminant $(a-a_i)^2+4c\bar c$ is positive we can distinguish $\alpha_i$ and $\beta_i$. Then one takes  $x(u)=A(u-\alpha_1)\cdots (u-\alpha_{k+1})$ and $y(u)=B(u-\beta_1)\cdots (u-\beta_{k+1})$
with $AB=c^{k+1}$ and $a,c$ and $A/\vert A\vert$ give 4 real local coordinates for the metric. 

The ${\bf D_k}$ case is more complicated -- here $z(u)$ is a quartic polynomial. Rewriting the equation as 
$$zx^2-(zy+\prod_{i=1}^ka_iu)^2=-\prod_{i=1}^k(z+a_i^2u^2)$$
we observe that a solution $x(u),y(u)$ means that 
\begin{equation}
f(u, w)=\frac{\sqrt{z(u)}x(u)-(z(u)y(u)+u^k\prod_{i=1}^k  a_i)}{\sqrt{z(u)}x(u)+(z(u)y(u)+u^k\prod_{i=1}^k  a_i)}
\label{function}
\end{equation}
is a rational function on the elliptic curve $w^2=z(u)$ with zeros $p_i$ and  poles $q_i$ at $w=\pm ia_ju$ for some choice of signs. This is a single constraint,  implemented by an elliptic integral,  on the divisor class $D=\sum_i p_i-q_i$ which reduces  the five coefficients of $z(u)$ to produce 4 local coordinates on the ALE space.

\begin{rmk} Since this approach predates Kronheimer's work it is perhaps useful to point out the context. Both cases produced formulae in terms of local coordinates, and in the ${\bf A_k}$ case one could relate this in \cite{Hit1} to a previously derived expression for the metric \cite{GH} by Gibbons and Hawking whose global interpretation was clear. In that sense the contribution of the twistor approach was to describe all the complex structures in the hyperk\"ahler family. At the time, there was no alternative expression for the ${\bf D_k}$ case and the local coordinates used were not sufficient to give a global description of the space, but much later it made sense when moduli spaces of monopoles and their hyperk\"ahler metrics were being analyzed. The constraint on the elliptic spectral curve of a charge 2 $SU(2)$ monopole with singularities was almost the same as the one above. The difference, which changes the asymptotic behaviour of the metric to ALF (asymptotically locally flat) is that the divisor class $D$, instead of being trivial and giving a rational function, is  the restriction of a  canonical element in $H^1(T\PP^1,\OO)$ to the curve $w^2=z(u)$ where $T\PP^1=\OO(2)$ is a non-compact  surface, the tangent  bundle of  $\PP^1$. The formulae in \cite{CH} require only slight modification to apply to the ALE case. 
\end{rmk} 

Proceeding to the ${\bf E_k}$ case following this approach proved difficult  and it becomes easier, but less explicit, to consider the algebraic geometry of compactifications.
\section{Compactifications}\label{compact}
If we compactify $\C^2$ to $\PP^2$ then $\C^2/\Gamma$ can be compactified to $\PP^2/\Gamma$ by adding the orbifold $\PP^1/\Gamma$ at infinity to the singular surface and its resolution. For the ${\bf {D,E}}$ series there are three orbifold points which give singularities of $\PP^1/\Gamma$ corresponding to the stabilizers of the vertices, midpoints of faces and edges of a  dihedron, tetrahedron, octahedron and icosahedron. These can be resolved   \cite{Pink} to give a configuration of rational curves which support an anticanonical divisor: $-K=2C+E+F+G$. For  ${\bf E_k}$  it is:

\hskip 4cm\begin{tikzpicture}
\draw[] (0,3) -- (6,3);
\draw[] (1,3) -- (1,0);
\draw[] (3,3) -- (3,0);
\draw[] (5,3) -- (5,0);
\node at (0.5,0) {${ E}$};
\node at (2.5,0) {${ F}$};
\node at (4.5,0) {${ G}$};
\node at (4.5,0) {${ G}$};
\node at (0.25,1.5) {${ 3}-{ k}$};
\node at (2.5,1.5) {$-{ 2}$};
\node at (2.0,3.5) {${C}$};
\node at (3.5,3.5) {${- 1}$};
\node at (4.5,1.5) {${ -3}$};
\end{tikzpicture}

and for ${\bf D_k}$ and ${\bf A_{2\ell-1}}$:

\begin{tikzpicture}
\draw[] (0,3) -- (6,3);
\draw[] (1,3) -- (1,0);
\draw[] (3,3) -- (3,0);
\draw[] (5,3) -- (5,0);
\node at (0.5,0) {${ E}$};
\node at (2.5,0) {${ F}$};
\node at (4.5,0) {${ G}$};
\node at (4.5,0) {${ G}$};
\node at (0.1,1.5) {$-{(k-2)}$};
\node at (2.5,1.5) {${ -2}$};
\node at (2.0,3.5) {${C}$};
\node at (3.5,3.5) {${ -1}$};
\node at (4.5,1.5) {${ -2}$};
\end{tikzpicture}
\hskip 2cm\begin{tikzpicture}
\draw[] (0,3) -- (4,3);
\draw[] (1,3) -- (1,0);
\draw[] (3,3) -- (3,0);
\node at (0.5,3.5) {${C}$};
\node at (0.5,1.5) {${- \ell}$};
\node at (2.7,1.5) {${- \ell}$};
\node at (2.3,3.5) {${ 0}$};
\node at (0.5,0) {${ D_1}$};
\node at (2.6,0) {${ D_2}$};
\end{tikzpicture}

\begin{figure}
\hskip 4cm\begin{tikzpicture}
\draw[] (4,3) ellipse (3cm and 1.5cm);
\draw[] (1,0) -- (1,6);
\node at (1,6.5) {${G}$};
\node at (0.75,3) {${ f}$};
\node at (0.75,0) {${ x}$};
\node at (4,1.25) {${ E}$};
\node at (1.75,4.35) {${ e}_1$};
\node at (3,4.7) {${ e}_2$};
\node at (4.25,4.75) {${ e}_3$};
\node at (5.7,4.5) {${ e}_4$};
\end{tikzpicture}
\caption{$D_4$ blowing up}
\label{pic}
\end{figure}
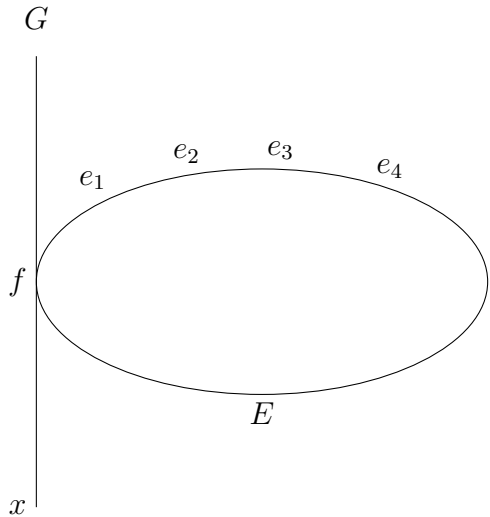

\begin{rmk}
We consider ${\bf A_k}$ with $k$ odd for uniformity here. The point is that the ${\bf D_k,E_k}$ Dynkin diagrams have a unique trivalent  vertex which means the corresponding rational curve in the resolution of the singularity is pointwise fixed. This will play an important role later. In the ${\bf A_k}$ case there is a further  circle action -- a circle subgroup of $SU(2)$  commuting with the action of the cyclic group $\Gamma$ and some combination fixes any one of the chain of rational curves. It is convenient to take the central vertex as the fixed one when $k$ is odd. In addition $-1\in \Gamma$ for all cases except ${\bf A_k}$ for $k$ even.
\end{rmk}
Conversely a rational surface with a canonical divisor as above is the compactification of one of the versal deformations of the singular surface. This way the surface and its subvarieties replaces the explicit equation. Following \cite{Hit3} consider the  configuration in $\PP^2$ in Figure \ref{pic}: a conic $E$ with 4 points $e_1,e_2,e_3,e_4$ and a  line through $x$ tangent to the conic at $f$. Blowing up the six points gives a cubic surface (in fact the natural projective compactification of the ${\bf D_4}$ equation (\ref{Dequation})). Because $E$ meets $G$ tangentially in the cubic we have three lines through a point and blowing up that point to give a $-1$ curve $C$ gives three $-2$ curves and the compactification above for ${\bf D_4}$. Clearly blowing up   further points $e_5, e_6, \dots, e_k$ gives the ${\bf D_k}$ configuration. 
The ${\bf E_k}$ series is similar but with an additional point $y$ on $G$ blown up. 

The above describes a compactification of the fibres of the twistor space $Z\rightarrow \PP^1$. The ALE property of the metric means that the 4-manifold has a conformal orbifold compactification by a single point (as used in \cite{At}) and in consequence its twistor space has an orbifold  compactification by adding $\PP^1/\Gamma$. After blowing up this gives a compactification  which extends the fibration so we obtain $\bar Z\rightarrow \PP^1$ where each fibre $\bar Z_u$ for $u\ne 0,\infty$ is holomorphically equivalent to one of the projective surfaces described above.

\section{The $\C^*$-action}
The circle action on the ALE space induces a corresponding action on the twistor space which extends to a $\C^*$-action, covering $u\mapsto \lambda u$  under the projection $\bar Z\rightarrow \PP^1$. It preserves the fibres over $u=0,\infty$ which correspond to complex structures $\pm I$, the resolution of the singularity and its complex conjugate. On the resolution 
the action fixes pointwise the rational curve defined by the trivalent vertex (we shall call this $C_0$) and the intersection points of the other rational curves. We shall consider the orbit of a generic twistor line $P$ under this action. 

Firstly consider flat space $\R^4$. The twistor space $Z$ here is the total space of the vector bundle $\OO(1)\oplus \OO(1)\rightarrow \PP^1$ and its compactification by $\PP^1$ is $\PP^3$. Blowing up this line at infinity gives $\bar Z\rightarrow \PP^1$ with fibres $\PP^2$.  In homogeneous coordinates we remove the line $z_0=z_1=0$ from $\PP^3$ and the projection to $\PP^1$ is $u=z_1/z_0$. The 4-parameter family of twistor lines considered as sections of $\OO(1)\oplus \OO(1)$ are
$z_2=az_1+bz_0, z_3=cz_1+dz_0$. The real lines are given by $d=\bar a, c=-\bar b$. Fix $(a,b,c,d)$ to give a twistor line $P$, then $\lambda\in \C^*$ transforms $P$ to
$$z_2=\lambda az_1+\lambda^{-1}bz_0,\qquad z_3=\lambda c z_1+\lambda^{-1}dz_0.$$
Eliminating $\lambda$ we see that the closure of the orbit of $P$ is the quadric surface in $\PP^3$ given by 
$$(cz_2-az_3)(dz_2-bz_3)+(ad-bc)z_0z_1=0.$$
In particular the surface meets the line at infinity in two points  $[0,0,a,c],[0,0,b,d]$.  
The quadric meets the fibre $u=0$ where $z_1=0$ which from the equation reads $(cz_2-az_3)(dz_2-bz_3)=0$ and consists of  two lines $L_1,L_2$ meeting at the origin, the closure of two orbits in $\C^2$ under the scalar $\C^*$-action. One line is given by the orbit of  the intersection of $P$ with the fibre $u=0$, the other is a line in $\PP^3$, the limiting transform  as $\lambda \rightarrow \infty$. 
 
 In the ALE case we take a generic twistor line $P$ and consider the closure $X$ of the orbit.  The orbifold compactification of $Z$ clearly produces a similar pair of lines $L_1,L_2$ as the flat case. Now consider the $\C^*$-action on the complement in $\bar Z$ of the fibre over $u=\infty$, as in \cite{Hit3} Section 7. The open Bialynicki-Birula stratum for the action consists of orbits which flow down to the fixed curve $C_0$ in the fibre over $u=0$. Thus a generic point on $P$ has a limit on $C_0$ and so this curve is another  component of the intersection $X\cap \bar Z_0$.   Note  that the intersection number $(L_1+L_2+C_0).C_0=1+1-2=0$. In fact  the cohomology class of $X$ in $H^2(\bar Z,\Z)$ must vanish on $C_0$ since the class $[C_0]\subset H_2(\bar Z,\Z)$ appears as a holomorphic curve with conjugate complex structures in $u=0,\infty$. 
 
 The curve $C_0$ is the limit for the downward flow of a generic point in $\bar Z_0$ and so the lines $L_1,L_2$, orbit closures, only intersect $C_0$ and not the other rational curves in the resolution. Then the cohomology class of $X$ vanishes on all classes in $H_2(Z_0,\Z)$.

 Now consider the intersection $X\cap \bar Z_1$, the fibre over $u=1$. Since $H_2(Z_1,\Z)\cong H_2(Z_0,\Z)$, we  learn from the above discussion that the divisor class $Q$ of $X$ in $\bar Z_1$ has the property that $Q$ is trivial on the affine part $H_2(Z_1,\Z)\subset H_2(\bar Z_1,\Z)$ and also on $E,F,G$ for the ${\bf D_k,E_k}$ series and $D_1,D_2$ for ${\bf A_{2\ell-1}}$. This yields the following:
\begin{thm} \label{Q} The divisor class $Q$ of the image of a twistor line in the rational surface $\bar Z_1$ is given by 
\begin{align}
{\bf E_6}: Q &= 12C+4E+6F+4G \nonumber \\
{\bf E_7}: Q &= 24C+6E+12F+8G \nonumber \\
{\bf E_8}: Q &= 60C+12E+30F+20G \label{divisor}\\
{\bf A_{2\ell-1}}: Q &= \ell C+D_1+D_2 \nonumber \\
{\bf D_k}: Q &= (2k-4)C+2E+(k-2)F+(k-2)G.\nonumber
\end{align}
and  in each case the  intersection number  $Q^2=\vert \Gamma\vert$ and $KQ=-4$ where the ALE space is asymptotic to $\C^2/\Gamma$.
\end{thm} 
\begin{prf} Since the curves $C,E,F,G$ and $D_1,D_2$ are disjoint from the representative cycles in $H_2(Z_0,\Z)\cong H_2(Z_1,\Z)$ we can look for classes $Q=a_0C+a_1E+a_2F+a_3G$ etc. which satisfy $QE=QF=QG=0$ and also, since $L_1,L_2$ intersect $C$ transversally, $QC=2$. Using the configurations in Section \ref{compact} this leads to the list. Then $Q^2$ and $KQ$ can be calculated directly using the adjunction formula for the rational curves $C,E,..$ etc.
\end{prf}
Since the image of a twistor line is a rational curve, the adjunction formula $KQ+Q^2=2g-2$ says that these curves are singular and if they have simple nodes then there must be $\vert \Gamma\vert/2 -1$ singularities. Thus each twistor line defines a nodal curve in the linear system $Q$. However, Honda and Nakata \cite{Honda1} prove that if a smooth projective  surface $S$ has a rational curve $C$ with $\delta>0$ nodes as all its singularities and $C^2+1-2\delta>0$, then the variety of such curves is non-singular and $C^2+1-2\delta$-dimensional. In our case $Q^2+1-2\delta=3$ and so the images of the 4-parameter family of twistor lines give an open set in the 3-parameter variety of nodal curves in the linear system.
The context of Honda's work appears in Section \ref{honda}. We shall henceforth denote the surface $\bar Z_1$ by $S$.
\begin{rmk} Note that since $E$ came from a conic in the plane and $G$ a line under the blowing up from the plane, on blowing down a nodal curve in the linear system $Q$ to $\PP^2$ the coefficients of $E,G$ in equations (\ref{divisor}) show that the image is a plane curve in general of high degree and with singular behaviour at $f$, apart from the movable nodes. 
\end{rmk}
\section{Examples} 
\subsection{ ${\bf A_{2\ell-1}}$ }\label{Akcurves}
The configuration of curves for ${\bf A_{2\ell-1}}$ can be achieved by taking $S$ to be the projective bundle $\PP(\OO\oplus \OO(\ell))\rightarrow \PP^1$ with $2\ell$ points $a_1,\dots, a_{2\ell}$ on the zero section blown up. Then $C$ is a fibre, $D_1$ the infinity section and $D_2$ the blown up zero section. 

\hskip 5cm\begin{tikzpicture}
\draw[] (0,3) -- (4,3);
\draw[] (1,3) -- (1,0);
\draw[] (3,3) -- (3,0);
\node at (0.5,3.5) {${C}$};
\node at (0.5,1.5) {${- \ell}$};
\node at (2.7,1.5) {${- \ell}$};
\node at (2.3,3.5) {${ 0}$};
\node at (0.5,0) {${ D_1}$};
\node at (2.6,0) {${ D_2}$};
\end{tikzpicture}

From the equation (\ref{divisor}) we have $Q = \ell C+D_1+D_2$ and $QC=2$ which means that, if $x$ is an affine coordinate on $\PP^1$,  $Q$ has an equation of the hyperelliptic form 
$$w^2+p(x)w+q(x)=0$$ where $p(x)$ has degree $\ell$ and $q(x)$ degree $2\ell$ and vanishes at each $x=a_i$, so $q(x)=c(x-a_1)\dots (x-a_{2\ell})$. 
Here  $\Gamma$ is the cyclic group $\Z_{2\ell}$, so we require $\vert \Gamma \vert/2 -1=\ell-1$ nodes. This means $p(x)^2-4q(x)$ must have $\ell-1$ double zeros. This is $\ell-1$ constraints on the $\ell+1$ coefficients of $p(x)$ together with  the multiplier $c$ giving the expected $(\ell+2)-(\ell-1)=3$ parameters.

Note that for $\ell=1$ there are no singularities and the image of a twistor line  is a  smooth rational curve of self-intersection $2$. In the model here we can blow down the $-1$ curves $D_1$ and $D_2$ to obtain a quadric surface and then the curves are plane sections. We shall meet this fact again in the final section.

When $\ell=2$, by a projective transformation  we can take three of the four points to be $(0,1,\infty)$ and write $q(x)=\lambda^2 \wp(z)$ using the Weierstrass $\wp$-function. Then the general constraint is for $\lambda \wp'(z)+a_0+a_1\wp(z)+a_2\wp^2(z)$ to have a double zero. 
\subsection{ ${\bf D_4}$ }
In this case equation (\ref{divisor})  gives $Q=4C+2E+2F+2G=-2K$. Now, as discussed in Section  \ref{compact},  $S$ is a cubic surface $\Sigma\subset \PP^3$ blown up at a point $a$ giving the $-1$ curve $C$. Since the embedding of $\Sigma$ is the anticanonical system, curves in the linear system $Q$  consist of intersections with quadric surfaces in $\PP^3$ meeting the cubic tangentially at that point and three others. We shall see next in this approach the constrained elliptic curve in the description of Section \ref{twistor}.

The cubic surface $\Sigma$ is the blow-up of the plane at the six points $x,f,a_1,a_2,a_3,a_4$ as in the diagram. We denote the images in $\Sigma$ of the curves $E,F$ etc. in $S$ by $E_*,F_*$ etc. and a generic twistor line $P\cong \PP^1$ maps to the nodal curve $Q_*\subset \Sigma$. The lines $E_*,F_*,G_*$ meet at the point $a$. Blowing down further to the plane we write the curves as $E_{**}$ etc.
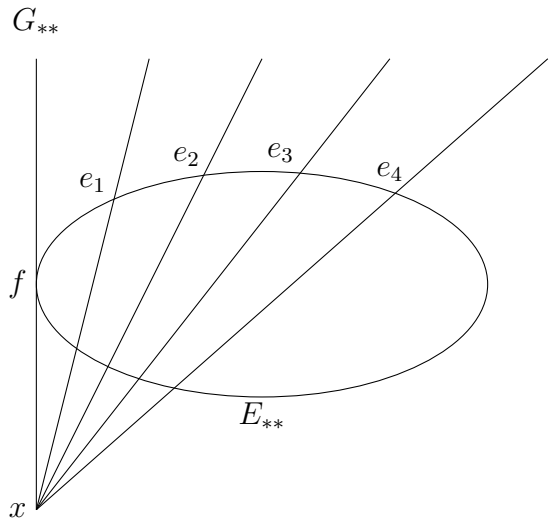
\begin{figure}
\hskip 4cm\begin{tikzpicture}
\draw[] (4,3) ellipse (3cm and 1.5cm);
\draw[] (1,0) -- (1,6);
\draw[] (1,0) -- (2.5,6);
\draw[] (1,0) -- (4,6);
\draw[] (1,0) -- (5.7,6);
\draw[] (1,0) -- (7.8,6);
\node at (1,6.5) {${G_{**}}$};
\node at (0.75,3) {${ f}$};
\node at (0.75,0) {${ x}$};
\node at (4,1.25) {${ E_{**}}$};
\node at (1.75,4.35) {${ e}_1$};
\node at (3,4.7) {${ e}_2$};
\node at (4.25,4.75) {${ e}_3$};
\node at (5.7,4.5) {${ e}_4$};
\end{tikzpicture}
\caption{the planar model}
\label{pic2}
\end{figure}
Consider the pencil of planes in $\PP^3$ containing the line $E_*\subset \Sigma$. This is the linear system $-K_{\Sigma}-E_*$, and 
$$Q_*(-K_{\Sigma}-E_*)=(-2K)(-K+C-E)=4.$$
A plane through $E_*$ meets the line $X_*$ in a point so $X_*$ parametrizes the pencil and $Q_*$ is a 4-fold cover of $X_*$. 
 
In  Figure \ref{pic2}, the lines through $x$ are the points of $X_*$, but also express the conic $E_{**}$ as a double cover of $X_*$  branched over $f$ and the point of contact of the other tangent  from $x$. Let $t$ be an affine parameter on $E_*$ where $t=0$ is the point $f$ and the covering involution is $t\mapsto -t$. Then $t^2$ is an affine parameter on $X^*$ and $t^2=0$ is the  line $G_{**}$. The plane in the pencil defined by $t^2=0$ meets $Q_*$ in the image of the zeros of a quartic polynomial $z(u)$ on $\PP^1$ and then, if $t=a_i$ is the point $e_i$ on the conic $E_{**}$, $t=-a_i$ is the other intersection of the line joining $x$ to $e_i$. The rational function (\ref{function}) on the elliptic curve $t^2=z(u)$ is then
$$f(t)=\prod_{i=1}^4\frac{(t-a_i)}{(t+a_i)}.$$
\section{The central sphere}
We take a small diversion here to see a rather different interpretation of the surface $S=\bar Z_1$. The element $-1\in S^1$ defines a holomorphic involution on the twistor space and composing with $\sigma$, a new real structure $\tau$ which covers $u\mapsto 1/\bar u$ with fixed point set $\vert u\vert=1$ -- the reflective rather than antipodal real structure on $\PP^1$. So the fibre $S$ given by $u=1$ has a real structure. The holomorphic curve $C_0\subset\bar Z_0$ is fixed by $-1$ and becomes a component of the real points of $S$. From this point of view $S$ is the complexification of 
this central 2-sphere and the paper \cite{Hit3} gave a description of the metric induced  from the hyperk\"ahler metric on the ALE space. 

The 2-sphere lies in an affine surface in $\R^3$ and inherits the natural volume form. Its conformal structure is defined by two complex conjugate foliations by null curves of the compexification. In \cite{Hit3} these are shown to be given by two conjugate pencils of rational curves of zero self-intersection. Moreover these curves are characterized by intersecting $C_0$ and $C$ in one point and having zero intersection with the cycles coming from the other rational curves in the resolution. Given our derivation of $Q$ in the previous section, it is clear that $Q=D+\tau(D)$ where $D$ is the divisor class for one pencil. 

\begin{ex} For ${\bf A_1}$, $S$ blows down to  a quadric surface $\PP^1\times \PP^1$ and the null curves are the two families of lines $D\cong \OO(1,0),\tau(D)\cong \OO(0,1)$. The projected twistor lines are plane sections, in the divisor class $\OO(1,1)$.  The Eguchi-Hanson metric on the cotangent bundle of the 2-sphere  is rotation-invariant and therefore a  standard metric on the real quadric $x_1^2+x_2^2+x_3^2=1$ in $\R^3$.
\end{ex} 

\begin{rmk}
An expression for $D$ is given in \cite{Hit3} in terms of the exceptional curves arising from blowing up the points $e_i,f,g,x,y$ on $\PP^2$ but there is an unfortunate choice of notation with $\bar E_i$ denoting the $-1$ curve coming from the line in the plane joining $x$ to $e_i$. If $E_i$ is the blow up of the point $e_i$ then in the ${\bf D_k}$ case $\bar E_i$ is indeed its conjugate but this does not necessarily hold for the ${\bf E_k}$ series. 

One way to determine the action of the real structure $\tau$ on the cohomology of $S$ is to note that, because it changes orientation on the holomorphic 2-spheres in the resolution, the action is $-1$ on the cohomology of the affine part $Z_1$. The circle action acts trivially on the $E,F,G$ configuration at infinity and so $\tau$ coincides with the antipodal real form. Then  the total class $E+F+G$ in $H^2(S,\Z)$ is acted on by $+1$ but in the case of ${\bf E_6}$ consideration of the real structure $(z_1,z_2)\mapsto (\bar z_1,-\bar z_1)$ acting on the invariant polynomials 
$$z_1^4+2i\sqrt{3}z_1^2z_2^2+z_2^4,\quad z_1^4-2i\sqrt{3}z_1^2z_2^2+z_2^4,\quad  z_1z_2(z_1^4-z_2^4)$$
shows that $\tau$ interchanges the cohomology classes $E$ and $G$.
\end{rmk}
\section{Einstein-Weyl manifolds}\label{honda}
We can now place the theorem we used of Honda and Nakata  in the context they were studying, which is to regard the complex surface $S$ as a minitwistor space for a 3-dimensional Einstein-Weyl manifold $M$. This is a manifold $M$ with a conformal structure and a torsion-free connection preserving it such that the symmetrized Ricci tensor $R_{(i,j)}$ is a scalar. The skew-symmetric component $R_{[i,j]}$ is a 2-form which is essentially the curvature of the line bundle of volume forms. This was the original idea in 4 dimensions of incorporating the electromagnetic field into general relativity. 
In practical terms one can choose a metric $g$ in the conformal equivalence class and then the Weyl connection satisfies $\nabla g=2\omega\otimes g$ for a 1-form $\omega$.

The differential-geometric link with the ALE spaces is the fact that the quotient of a hyperk\"ahler 4-manifold by an isometric circle action is an Einstein-Weyl 3-manifold. In local coordinates a 
representative metric can  \cite{Calder1} be put in  the form
\begin{equation} g=e^u(dx^2+dy^2)+dt^2,\qquad \omega=-u_tdt
\label{ew}
\end{equation}
where $u$ satisfies the $SU(\infty)$ Toda equation $u_{xx}+u_{yy}+(e^u)_{tt}=0$. 

We can describe this more globally, following \cite{Hit4}. The circle action on the ALE space preserves the K\"ahler form $\omega_1$ for this complex structure. Then $t$ in the equation is, up to a constant, the moment map for this action. The quotient by the circle action of a level set $t=c$ is a K\"ahler quotient, with an induced holomorphic structure  and $u_te^u(dx^2+dy^2)$ at $t=c$ is this metric. The time-dependent metric $h=e^u(dx^2+dy^2)$ on this surface satisfies the second order equation  $h_{tt}=\kappa h$ where $\kappa$ is the Gaussian curvature. 

On the hyperk\"ahler manifold the circle action generates a vector field $X$ and the vector field $IX$ is horizontal and is the gradient of the moment map  $t$. Thus the $dt^2$ term in formula (\ref{ew}) shows that this choice of metric in the conformal class is determined by rescaling the quotient metric with the norm square of $X$ to make $dt$ have length $1$. 

The {\emph {twistor theory }} for Einstein-Weyl 3-manifolds was established in \cite{Hit2} and based on the minitwistor quotient of a twistor space: if there is a good quotient of the complex 3-manifold $Z$ by a $\C$-action, then the normal bundle of a twistor line is an extension $\OO\rightarrow \OO(1)\oplus \OO(1)\rightarrow \OO(2)$ where the trivial subbundle is tangent to an orbit of the action. Then the quotient is a surface $S$ and the twistor lines project to  rational curves of self-intersection $2$. The two models for this are $S=\OO(2)$ the tangent bundle of $\PP^1$, and $S=\PP^1\times \PP^1$ with the corresponding Einstein-Weyl structures being $\R^3$ with the flat metric or hyperbolic metric respectively. Furthermore, the only compact complex surfaces which admit such curves give locally the same geometry. 

This is where the work of \cite{Honda1} comes to the rescue, considering families of nodal curves. In that paper a local tubular neighbourhood of a rational nodal curve is constructed which can be ``normalized" to the neighbourhood of a smooth curve of self-intersection $2$. Our surface $S$ is a global substitute for the quotient of the ALE twistor space -- locally around a twistor line we can find a quotient, and this will constitute  an example of this local normalization. The authors of  \cite{Honda1}, \cite{Honda2} apply this to  introduce many examples of  Einstein-Weyl manifolds from the geometry of compact rational surfaces. In fact, as the first author of \cite{Honda1} has pointed out,  Section 5.2 of  that paper contains a family which includes our ${\bf A_{2\ell-1}}$ case. Note that there are now two real structures $\sigma$, $\tau$ on the surface $S$. The antipodal structure $\sigma$ gives an Einstein-Weyl conformal structure which is positive-definite whereas $\tau$, whose fixed points give the real surface $C_0\subset S$ for the ALE metric, gives the quotient Lorentzian signature.

In the case of    ${\bf A_{2\ell-1}}$  
the surface $S$ and its nodal curves have already been described here in Section \ref{Akcurves}. There is another action -- a   triholomorphic action  -- of the circle on the 4-manifold, and the  quotient   Einstein-Weyl structure is here just the Levi-Civita connection of the flat metric on $\R^3$.  This is essentially the Gibbons-Hawking Ansatz \cite{GH}.
We can see this from the twistor point of view: the $\C^*$-action $(x,y,z)\mapsto (\lambda x,\lambda^{-1}y, z)$ on $\OO(k+1)
\oplus \OO(k+1)\oplus \OO(2)$ preserves $z$, and twistor lines project to  sections $z(u)$ of $\OO(2)\rightarrow \PP^1$, the minitwistor space for $\R^3$. 

For the action we are considering in this paper an explicit formula for the ${\bf A_k}$ Einstein-Weyl structure can be found in \cite{Calder1}, based on earlier work of Ward \cite{Ward}.  For $k=1$, we already observed that $S$ is a quadric and so the Einstein-Weyl quotient is hyperbolic 3-space, but it is useful to see this concretely.  One standard form of this ${\bf A_1}$ metric is the following, using left-invariant 1-forms $\sigma_i$ on $SU(2)$:
$$g=\frac{dr^2}{(1-(1/r)^4)}+\frac{r^2}{4}(\sigma_1^2+\sigma_2^2)+\frac{r^2}{4}(1-(1/r)^4)\sigma_3^2$$
which clearly shows that as $r\rightarrow \infty$ it becomes asymptotically a quotient of Euclidean space and as $r\rightarrow 1$ we get the 2-sphere $C_0$ with metric $(\sigma_1^2+\sigma_2^2)/4$. The vector field $X$ given by the circle action is dual to $\sigma_3$.

If we take the K\"ahler form to be  
$2rdr\wedge\sigma_3+
{r^2}\sigma_1\wedge \sigma_2$
then  the moment map is $t=r^2$ and  $(X,X)=(t-t^{-1})/4$. This gives the data for the Toda form for the metric  (\ref{ew}) and one finds $\omega=-d\log(t^2-1)$ and so $d\omega=0$. This means the volume bundle for the conformal structure is flat and we can rescale the metric by a power of $(t^2-1)$ to obtain a Levi-Civita connection for  the hyperbolic 3-space:   
$$\frac{\rho^2}{(1-\rho^4)^2}\left(d\rho^2+\frac{1}{4}\rho^2(\sigma_1^2+\sigma_2^2)\right).$$
Here $2r^2=\rho^2+\rho^{-2}$ and so $\rho\rightarrow 1$ corresponds to $r\rightarrow 1$ -- the central 2-sphere becomes the 2-sphere at infinity in hyperbolic space. As $r\rightarrow \infty$, $\rho\rightarrow \infty$ and infinity in the ALE space defines a distinguished origin in the hyperbolic space. 

These two features hold in the general ALE case: the Einstein-Weyl quotient has an asymptotic 2-sphere and an orbifold central point modelled on $\R^3/\Gamma$ where $\Gamma$ is acting via its image $\Gamma/\pm 1$ in $SO(3)$. There will be other orbifold singularities pertaining to isolated fixed points of the circle action. This is not the place to go into further details but it allows for the opportunity to give a geometrical meaning to the appearance of nodes in the twistor lines. 

In the original situation of smooth rational curves of self-intersection $2$, the space $S$ is interpreted as the space of oriented geodesics for the Weyl connection, the antipodal real structure being the change of orientation. Then two real twistor lines intersect in two points because two points in the 3-manifold have a unique geodesic connecting them (with two orientations). 
In our case the self-intersection $Q^2=\vert \Gamma\vert$ so we expect $\vert\Gamma\vert/2$ geodesics between them. In a neighbourhood of the orbifold point this is clear: take a fundamental domain in the $\Gamma/\pm 1$ - covering and two  points $p,q$. Then we shall find geodesics joining $p$ to each of the $\Gamma$-translates of $q$. 

The same principal applies to a single point $p$. In a local normalized neighbourhood the smooth rational curve represents the geodesics corresponding to the tangent directions at $p$ along each of which passes a geodesic. But, as with the geodesics on a cone, globally we obtain geodesics which return with a different direction, giving an identification of pairs of points on the curve, and hence a node on its image in the space of geodesics. We lose one of the geodesics when we had two points to give $\vert \Gamma\vert/2-1$ singularities.

\vskip 1cm
 Mathematical Institute,  Woodstock Road, Oxford OX2 6GG, UK
 
 hitchin@maths.ox.ac.uk


\begin{thebibliography} {11} 
 \bibitem{At} 
M.F.Atiyah, {\it Green's functions for self-dual four-manifolds},  Adv. Math. Suppl. Stud. {\bf 7} (1981) 129 --158.

\bibitem{Calder1}
D.Calderbank, {\it The geometry of the Toda equation}, J. Geometry and Physics {\bf 36} (2000) 152--162.



\bibitem{CH}
S.Cherkis \& N.Hitchin, {\it Gravitational instantons of type $D_k$}, 
Comm.Math.Phys. {\bf 260} (2005) 299--317. 


\bibitem{Pink}
M.Demazure, H.Pinkham \& B.Teissier (eds.), ``S\'eminaire sur les singularit\'es des surfaces", Lecture Notes in Mathematics {\bf 777}, Springer Verlag Heidelberg (1980).
%

\bibitem{GH}
G.Gibbons \& S.Hawking, {\it Gravitational multi-instantons}, 
Phys.Lett. {\bf B 78} (1978) 430--432. 


\bibitem{Hit1}
N.J.Hitchin,  {\it Polygons and gravitons}, 
Math.Proc.Camb.Phil.Soc.  {\bf 85} (1979) 465--476. 

\bibitem{Hit2}
N.J.Hitchin, {\it Complex manifolds and Einstein's equations},
in ``Twistor geometry and nonlinear systems'' , H.Doebner et al (Eds), 
Lecture Notes in Mathematics {\bf 970}, 73-99, Springer, Berlin (1982).

\bibitem{Hit4}
N.J.Hitchin, {\it Higgs bundles and diffeomorphism groups}, in ``Surveys in Differential Geometry  Vol 21" H-D Cao and S-T Yau, (eds.),  
International Press, Cambridge, Mass. (2016), 139--163

\bibitem{Hit3}
N.J.Hitchin, {\it The central sphere of an ALE space}, Quarterly Journal of Mathematics {\bf 72}  (2021) 277--337.


\bibitem{Honda1}
N.Honda and F.Nakata, {\it Minitwistor spaces, Severi varieties, and Einstein--Weyl structure}, Ann. Glob. Anal. Geom. {\bf 39} (2011) 293--323.

\bibitem{Honda2}
N.Honda, {\it Segre quartic surfaces and minitwistor spaces}, New York J. Math.  {\bf 28} (2022) 672--704.

\bibitem{Kron1}
P.B.Kronheimer, {\it The construction of ALE spaces as hyperk\"ahler quotients}, J.Differential Geometry {\bf 29} (1989) 665--683.

 \bibitem{Kron2}
P.B.Kronheimer, {\it A Torelli-type theorem for gravitational instantons}, 
J.Differential Geom. {\bf  29} (1989) 685--697. 


  \bibitem{Ward}
  R.S.Ward, {\it Einstein-Weyl spaces and $SU(\infty)$ Toda fields}, Class. Quantum Grav. {\bf 7} (1990) 95--98.
\end{thebibliography}
\end{document}